\documentclass{amsart}
\usepackage{changebar}
\usepackage{lineno}
\usepackage{amsthm}
\usepackage{verbatim}
\usepackage[utf8]{inputenc}
\usepackage{blindtext}

\makeatletter
\renewcommand\part{%
   \if@noskipsec \leavevmode \fi
   \par
   \addvspace{4ex}%
   \@afterindentfalse
   \secdef\@part\@spart}

\def\@part[#1]#2{%
    \ifnum \c@secnumdepth >\m@ne
      \refstepcounter{part}%
      \addcontentsline{toc}{part}{\thepart\hspace{1em}#1}%
    \else
      \addcontentsline{toc}{part}{#1}%
    \fi
    {\parindent \z@ \raggedright
     \interlinepenalty \@M
     \normalfont
     \ifnum \c@secnumdepth >\m@ne
       \Large\bfseries \partname\nobreakspace\thepart
       \par\nobreak
     \fi
     \huge \bfseries #2%
     \par}%
    \nobreak
    \vskip 3ex
    \@afterheading}
\def\@spart#1{%
    {\parindent \z@ \raggedright
     \interlinepenalty \@M
     \normalfont
     \huge \bfseries #1\par}%
     \nobreak
     \vskip 3ex
     \@afterheading}
\makeatother

\newtheorem{lemma}{Lemma}
\newtheorem{theorem}{Theorem}
\newtheorem{proposition}{Proposition}
\newtheorem{corollary}{Corollary}
\newtheorem{definition}{Definition}

\newtheorem{assumption}{Assumption}
\newtheorem{remark}{Remark}
\newtheorem{fact}{Fact}

\newcommand{\nth}[1]{$#1 {\rm - th }$}

\newcommand{\Z}{\mathbb{Z}}

\newcommand{\ZMs}[1]{(\Z / #1 \cdot \Z)^*}

\newcommand{\Tr}{\mbox{\bf Tr}}

\newcommand{\rg}[1]{\mbox{\bf #1}}
\newcommand{\eu}[1]{\mathfrak{#1}}
\newcommand{\id}[1]{\mathcal{#1}}
\newcommand{\Gal}{\mbox{ Gal }}

\newcommand{\Ker}{\mbox{ Ker }}

\newcommand{\rf}[1]{(\ref{#1})}

\newcommand{\nmid}{\not \hspace{0.25em} \mid}
\newcommand{\Norm}{\mbox{\bf N}}

\newcommand{\F}{\mathbb{F}}

\newcommand{\K}{\mathbb{K}}

\newcommand{\KL}{\mathbb{L}}
\newcommand{\M}{\mathbb{M}}
\newcommand{\Q}{\mathbb{Q}}

\newcommand{\C}{\mathbb{C}}
\newcommand{\E}{\mathbb{E}}
\newcommand{\N}{\mathbb{N}}

\def\ra{\rightarrow}
\newcommand{\ran}{\rangle}
\newcommand{\lan}{\langle}

\newcommand{\ol}{\overline}
\newcommand{\veps}{\varepsilon}

\newcommand{\st}{^{\times}}

\begin{document}

{\obeylines \small
\vspace*{-0.2cm}
 \hspace*{3.5cm} Par les quatre horizons qui crucifient le monde
 \hspace*{3.5cm} Par les quatre horizons qui crucifient le monde
 \hspace*{3.5cm} Par ceux qui sont sans pieds, par ceux qui sont sans mains
 \hspace*{3.5cm} Et par le juste mis au rang des assassins
 \hspace*{5.5cm} Je vous salue, Marie
\vspace*{0.2cm} 
\footnote{Francis James: {\em Pri\'ere}}
\vspace*{0.2cm} 

\vspace*{0.4cm}
\hspace*{4.5cm}\indent    {\it To all brothers and sisters, watching in awe}
\hspace*{5.5cm}\indent    {\it to things happening around}
\vspace*{0.5cm}
\smallskip}
\title[Fermat-Catalan] {The Strong Fermat-Catalan Equation} 
\author{Preda Mih\u{a}ilescu} 
\address[P. Mih\u{a}ilescu]{Mathematisches Institut der Universit\"at 
G\"ottingen}
\email[P. Mih\u{a}ilescu]{Preda@uni-math.gwdg.de} 
\date{Version 1.0 \today}
\vspace{0.3cm}
\begin{abstract}
We give a cyclotomic proof of the fact that the equation $\frac{x^p + y^p}{x+y}  =  p^e z^q$ has no 
solutions in coprime integers $x,y,z$ and $p > 3; q$, a pair of distinct odd primes. 
\end{abstract} 
\maketitle
\vspace*{-0.5cm}
\tableofcontents
\vspace*{1.0cm}
\section{Introduction}
The Fermat-Catalan equation
\begin{eqnarray}
    \label{fct}
    x^p + y^p = z^q; \quad x, y, z \in \Z; \quad (x, y, z ) = 1, \quad p \neq  q > 2
\end{eqnarray}
has been intensively investigated in the decades since Wiles's proof of Fermat's Last Theorem.
This paper does not aim to provide any overview of the existing results, and we refer to
\cite{BMS} for more details and literature. The methods used so far only succeeded
to prove that no solutions exist for a small set of primes $p$. 

We also consider the more general {\em Strong Fermat-Catalan} equation. 
\begin{eqnarray}
\label{sfct}
\quad \frac{x^p + y^p}{x+y} & = & p^e z^q; \quad (x, y, z) \in (\Z^*)^3, \quad ( x,y,z ) = 1, \\ \nonumber & &  
\hbox{and $p > 3; q$ are distinct primes; $e \in \{0, 1\}$ }.
\end{eqnarray}
The relation between the value of $e$ and $x, y, z$ is explained below. 
The connection to \rf{fct} will become apparent below, from the classical formulas of Barlow and Abel. 
We prove
\begin{theorem}
    \label{main}
    The Diophantine equation \rf{sfct} has no integer solutions. In particular, the Fermat-Catalan equation \rf{fct}
    has no solutions with distinct odd primes $p, q$ and $p > 3$. 
\end{theorem}
\begin{remark}
    \label{terms}
    We choose the term {\em Fermat-Catalan} for equation \rf{fct}, because of its sharing the 
    shape of a cyclotomic norm equation with the classical Catalan equation, and also involving
    two prime exponents. The terminology in this area is not well established, and some authors
    refer to the equation $x^p + y^q + z^r = 0$, involving three odd prime exponents, as
    Fermat-Catalan. For the reasons above, we consider that the term better applies to \rf{fct}.
    The name of {\em Strong Fermat-Catalan Eqution} echoes a designation introduced by Gras and Qu\^eme
    for the equations that relate to Fermat's equation in the same way as \rf{sfct} relates to\rf{fct}.
    To the best of our knowledge, this equation has not been considered separately in the literature to 
    this day, except for its particular case $y = -1$, in which one retrieves the classical equation of
    Nagell-Ljunggren.
\end{remark}

We note that for $p=3$, \rf{sfct} 
has infinitely many solutions that stem from norms of \nth{q} powers in the Eisenstein integers: 
if $(x, y) \in \Z^2$, $\rho$ is a complex third root of unity and 
$\xi^q = x + y \rho$ for some $\xi \in \Z[\rho]$, then $\frac{x^3+y^3}{x+y} = z^q$, 
with $z = \xi \cdot \bar{\xi} \in \Z^*$. The condition $p > 3$ is thus necessary. 

The proof below is almost identical to the proof given in a separate paper for the 
Strong Fermat Equation. We separated the cases $p = q$ and $p \neq q$ for ease of
the overview  and avoiding overloaded case distinctions. 
The reader having worked through the previous paper, should not be
surprised to find, not only analogous proofs, but entire sections of preparatory
results which coincide. 

\section{Classical results, notations and prerequisites}
A classical fact, often attributed to Euler, states that 
for coprime integers $x, y$ and $n \in 2\N+1,\quad n>1$, the greatest common divisor
$d = \left( \frac{x^n+y^n}{x+y}, x+y \right)$ divides $n$. In general, the
following holds:
\begin{fact}
    \label{euler}
    Let $\rg{K}$ be an abelian number field, let $x, y \in \id{O}(\rg{K})$ be coprime
    and $p$ be a rational prime. Then
    \begin{eqnarray}
        \label{eulergcd}
        \left(\frac{x^p + y^p}{x+y}, x+y\right) \ \big \vert \ p. 
    \end{eqnarray}
\end{fact}
\begin{proof}
This can be verified by using the substitution $s = x+y$ and 
introducing it in the fraction 
\[ \frac{x^p+y^p}{x+y} = \frac{x^p + (s-x)^p}{s} = \frac{p x^{p-1} + O( s )}{s} .\] 
Since $(x, s ) = 1$, the claim follows.
\end{proof}

In our case, if $d = 1$, then 
$p \nmid z$, while for $d = p$, the same substitution implies
that $v_p\left(\frac{x^p + y^p}{x+y}\right) = 1$, hence the introduction of $p^e$ in \rf{sfct}.

\begin{remark}
    \label{xpy}
    Since $( x+y, \frac{x^p+y^p}{x+y}) = 1$, it follows that $1 = (z^p, x+ y ) = (z, x+y)$, so $z$ and $x y( x+y)$ are coprime.
\end{remark}

Using the methods used for the proof of the Barlow and Abel relations -- \cite{Ri}, Lecture IV, \S 1 --
for the Fermat equations, we find the following factorizations with respect to \rf{fct}:
\begin{eqnarray}
\label{factf}
x + y  =  u^q &\hbox{and} & \quad \frac{x^p + y^p}{x+y} = v^q \quad \hbox{if $p \nmid (x+y)$, and} \\
\nonumber x+y  =  p^{q-1} \cdot u^q &\hbox{and} & \quad \frac{x^p + y^p}{x+y} = p v^q \quad \hbox{otherwise}.
\end{eqnarray}
We see that having solutions to the Strong Fermat-Catalan equation 
\rf{sfct} is a necessary -- 
but not sufficient -- 
condition for the Fermat-Catalan equation \rf{fct} to have solutions. Thus, results on solutions of \rf{sfct}
imply the same conditions for solutions of \rf{fct}. However, the case $p = 3$ of \rf{fct} requires 
a separate treatment, for reasons explained above.

In our proof, we make intensive use of a large submodule $J$ of the Stickelberger ideal $I$ introduced below, which 
comes with pleasant additional properties with respect to possible solutions of \rf{sfct}. Explicitly, for
$t \in J$ and $(x, y, z)$ a solution of \rf{sfct}, there is a $\beta( t ) \in \Z[ \zeta_p ]$ with
\[ \beta( t )^q = \frac{(y + \zeta_p x )^t}{p^{ek}}, \quad k \geq 1.\]
The numerical discriminant of $\beta( t )$ induces an intricate factorization of the invariant $K = x y (x+y)$, which is coprime to $z$, by Remark \ref{xpy} below.

\subsection{Notations}
Throughout this paper, for $r$ a prime or a prime power, we denote
by $\F_r$ the field with $r$ elements. 

We let $P = \{0, 1, \ldots, p-1\}, P^* = P \setminus \{ 0 \}$ be the minimal positive representatives 
for $\F_p$ and $\F_p^{\times}$, respectively; $\zeta$ will be a primitive \nth{p} root of unity and 
$\K = \Q[ \zeta ]$ the \nth{p} cyclotomic field, with Galois group $G = \Gal( \K /\Q )$. The automorphisms 
$\sigma_c \in G$ are given by $\zeta \mapsto \zeta^c$, for $c \in P^*$. The Teichm\"uller character is 
 \[ \varpi : G \ra \Z_p\st ; \quad \sigma_c \mapsto \varprojlim_n (c^{p^n} \bmod p^{n+1}) . \]
The character restricts to characters mapping $G$ to $\F_p$
and it can also be interpreted as a character of $\F_p\st$ via the isomorphism $G \cong \F_p\st$.
It induces the spectral decomposition of $\rg{R}[ G ]$, for $\rg{R} \in \{ \F_p, \Z_p \}$
 as follows: 
 \begin{eqnarray}
 \label{idemps}
 \veps_k & := & \frac{1}{p-1} \sum_{c \in P^*} \varpi^k( c ) \sigma_c^{-1}; \ k \in P^*, \quad \hbox{verifying:} \nonumber \\
 \sum_{k \in P^*} \veps_k & = & 1; \quad \quad \veps_k \cdot ( \sigma_c - \varpi( c )^k ) = 0.
 \end{eqnarray}
  The complex conjugation acting on $\K$ is denoted by 
 $\jmath = \sigma_{p-1} = \sigma^{(p-1)/2}$. 
 We use the uniformizor $\lambda = 1-\zeta \in \Z[ \zeta ]$,  that generates the principal prime 
 $\wp \subset \Z[ \zeta ]$ above $p$ (we remind that $\wp^{p-1} = (\lambda)^{p-1} = (p)$). 

We shall also fix a primitive \nth{q} root of unity $\xi \in \C$ and let $\K' = \Q[ \xi ]$ and 
\[ H = \Gal( \K'/\Q ) = \{ \tau_d \ : \ \xi \mapsto \xi^d; \ d = 1, 2, \ldots, q-1 \} \]
The composite field is $\KL = \Q[ \zeta, \xi ]$, an abelian extension for which we denote the canonical 
lifts of $G, H$ by $G', H'$ respectively, so $\Gal( \KL/\Q ) = G' \times H'$. Let $Q = \{ 0, 1, \ldots, q-1\}; \ Q^* = Q \setminus \{ 0 \}$ be defined with respect to $q$
and by analogy to the sets $P, P^*,$ etc. Let $\lambda' = 1-\xi$, so $\wp' = ( \lambda' )$
is the ramified prime of $\K'$ over $q$.

 \begin{assumption}
 \label{sol}
 We assume in the sequel that $(x, y, z) \in \Z^3$ is a solution to \rf{sfct}, for the prime exponents $p, q$, and $ x  > | y | \geq 1$. 
\end{assumption}

\subsection{Characteristic numbers and characteristic ideals}
In the cyclotomic field $\K$ we have:
\[ z^q = \frac{x^p + y^p}{p^{e}(x+y)}  = \prod_{c \in P^*} \frac{y + \zeta^c x}{(1-\zeta^c)^e}. \]
This leads naturally to the following 
\begin{definition}
\label{charids}
We define
\begin{eqnarray}
\label{charni}
 \alpha = \frac{y + \zeta x}{(1-\zeta)^e}; \quad \eu{A} = ( \alpha, z ) \subset \Z[ \zeta ],
\end{eqnarray}
as the {\em \textbf{characteristic number}} of the equation \rf{sfct} and $\eu{A} = ( \alpha, z )$ 
is {\em the \textbf{characteristic ideal}} of the equation.
\end{definition}
The Lemma \ref{aux1} below, shows that the characteristic number and ideal indeed encode the 
properties of solutions of \rf{sfct}.
\begin{lemma}
\label{aux1}
\begin{itemize}
\item[ 1. ] The characteristic number $\alpha$ is integral. 
\item[ 2. ] The Galois group $G$ acts on the characteristic number, giving raise to pairwise coprime integral elements; that is,  for $1 \leq c < d \leq p-1$,
\begin{eqnarray*}
( \sigma_c( \alpha ), \sigma_d( \alpha )) = (1) .
\end{eqnarray*}
\item[ 3. ] The {\em characteristic ideal} $\eu{A}$ as $ \eu{A} = ( \alpha, z )$ is related to the 
{\em characteristic number} $\alpha$ by the relations:
\begin{eqnarray}
\label{rclass}
\eu{A}^q & = & ( \alpha ), \quad \Norm( \eu{A} ) = (z).
\end{eqnarray}
\item[ 4. ] The characteristic number satisfies:
\begin{eqnarray}
\label{al'}
\frac{\alpha}{\bar{\alpha}} & = & \upsilon \cdot \frac{1 + \zeta ( x/y )}{1 + \bar{\zeta}(x/y)}, \quad \hbox{with } \\ \nonumber
\upsilon & = & \begin{cases}
      1 & \hbox{if $ e = 0$, and } \\
      -\bar{\zeta} & \hbox{for $e = 1$,}
    \end{cases}
\end{eqnarray}
and if $e = 0$ we have 
\begin{eqnarray}
\label{twist}
\alpha' := \zeta^{x/(x+y)} \alpha = (x+y)\cdot ( 1 + O(\lambda^2 )).
\end{eqnarray}
\end{itemize} 
\end{lemma}
\begin{proof}
Let's first prove point 1: for $e=0$, clearly $\alpha$ is an integral element, while if $e=1$, then $p | ( x+y)$, so $\alpha$ is also integral.

For point 2., let's first remind that $\lambda = 1 - \zeta$ and that, for distinct $c, d \in P$, we have $\frac{\zeta^c - \zeta^d}{\lambda}$ is a unit in $\Z[\zeta]$. Let $I( c,d ) = \left( \sigma_c( \alpha), \sigma_d( \alpha )\right) $; then $y \lambda \in I( c, d)$. If $e= 0$,
this follows from $\sigma_c( \alpha ) - \sigma_d( \alpha ) = ( \zeta^c -\zeta^d ) y \in I( c, d)$, and for $e= 1$,
we have $(1-\zeta^c) \sigma_c(\alpha) - (1-\zeta^d) \sigma_d(\alpha) = -( \zeta^c -\zeta^d ) y \in I( c, d)$. 
Likewise, $x \lambda \in I( c, d)$: for $e= 0$,
we have $\bar{\zeta}^c \sigma_c( \alpha ) - \bar{\zeta}^d \sigma_d( \alpha ) = ( \bar{\zeta}^c -\bar{\zeta}^d ) x \in I( c, d)$, 
while for $e= 1$,
we have $(1-\bar{\zeta}^c) \sigma_c( \alpha ) - (1-\bar{\zeta}^d) \sigma_d( \alpha ) =( \bar{\zeta}^d -\bar{\zeta}^c) x \in I( c, d)$. 
Recall that $\lambda = 1-\zeta$ and $\frac{\zeta^a - \zeta^b}{\lambda} \in \id{O}^{\times}( \K )$ for any distinct  $a, b \in P$. Thus
$I( c,d ) | ( x,y) ( \lambda ) = \wp$, since $(x,y) = 1$. However, $( \alpha, p) = (1)$ by definition, so it follows that $I( c, d ) = (1)$, as claimed.

For point 3, we multiply out the norm, to get $\Norm( \alpha ) = z^q$. So 
$\alpha | z^q = \prod_{c \in P^*} \sigma_c(\alpha)$, and thus 
$z^q/\alpha = \prod_{c \in P^*, \sigma_c \neq 1} \sigma_c(\alpha)$. 
Consequently $\left( \alpha, z^q/\alpha \right) = (1)$, by point 2. 
For the characteristic ideal, this implies: 
\[ \eu{A}^q = \left( \alpha^q, \alpha^{q-1} z, \ldots, \alpha z^{q-1}, 
\alpha \cdot (z^q/\alpha) \right) = ( \alpha ) \cdot J , \]
where the ideal $J = ( \alpha^{q-1}, \ldots, z^{q-1}, z^q/\alpha ) = ( \alpha, z^q/\alpha ) = (1)$, 
hence $\eu{A}^q = (\alpha)$: the characteristic
ideal is either principal or it has order $q$\footnote{The order of an ideal is naturally defined as the order of its class in the class group.}. 
The relation $\Norm( \eu{A} ) = ( z )$ follows from $\Norm( \alpha ) = z^q$. 

For point 4, \rf{al'} follows from $\frac{1-\bar{\zeta}}{1-\zeta} = - \bar{\zeta}$. For
\rf{twist} we note that for $c \in P^*$, 
\[ \zeta^c = (1 - \lambda )^c = 1 - c\lambda + O( \lambda^2 ), \]
and $y + \zeta x \equiv (x+y) - \lambda x \equiv (x+y) \cdot \zeta^{x/(x+y)} \bmod \lambda^2$. 

\end{proof}
\subsection{The Stickelberger ideal and its action}
The Stickelberger element $\vartheta = \frac{1}{p} \sum_{c=1}^{p-1} c \sigma_c^{-1} \in \frac{1}{p} \Z[ G ]$ generates
the Stickelberger ideal  in the group ring of $G$ over the rational integers, by intersecting its principal ideal with $\Z[G ]$, 
according to 
\begin{eqnarray}
\label{stickid}
 I = \vartheta \Z[ G ] \cap \Z[ G ]. 
 \end{eqnarray}
 
 Comparing to the definition of the orthogonal idempotents in \rf{idemps}, we note that
 \begin{eqnarray}
 \label{st-el}
 \vartheta &=& \frac{p-1}{p} ( \veps_1 - A p ), \quad A \in \Z_p[ G ].
 \end{eqnarray}
 \subsubsection{Generators and relations}
The ideal $I$ has the property of annihilating the class group of $\K$ ( \cite{Wa}, \S 15.1). 
That is, for each ideal $\eu{C} \subset \Z[ \zeta ]$ and each $\theta \in I$, 
the ideal $ \eu{C}^{\theta} \subset  \Z[ \zeta ]$ is principal. There exists a base for 
$I^- = (1-\jmath ) I$, made of $(p-1)/2$ elements, called {\em Fueter elements}, \cite{Fu}, see also \cite{Mi2}, 
which are
\begin{eqnarray}
\label{fueter} 
& & \psi_n =  \vartheta ( 1+\sigma_n -\sigma_{n+1}) = \sum_{c \in S_n } n_c \sigma^{-1}_c  \in \Z_{\geq 0}[ G ], \quad \hbox{for } n \in \left\{ 1, 2, \ldots, \frac{p-1}{2} \right\}\\ 
& & \hbox{with} \quad \nonumber n_c = \left( \left[ \frac{(n+1)c}{p} \right] - \left[ \frac{nc}{p} \right] \right) \\
& & \label{Fueter1or0} \hbox{and} \quad  n_c + n_{p-c} = 1,
\end{eqnarray}
where the support $S_n \subset \{1, 2 \ldots, p-1 \}$, satisfies\footnote{The set $p-S_n$ designates naturally $\{ p - r \ : \ r \in S_n \} $.}  
$S_n \cup ( p - S_n ) = P^*$ and
is deduced from the definition of $\psi_n$. We note that \rf{Fueter1or0} implies that, for $\psi_n = \sum_{c \in S_n } n_c \sigma^{-1}_c$ an element of the Fueter base, $n_c = 0$ or $n_c=1$, as well as $(1+\jmath) \cdot \psi_n = \Norm_{\Q(\zeta)/\Q}$. The multiples of the norm are thus the only 
elements of $(1+\jmath) I$. We shall denote conjugates $\psi = \sigma \psi_n$ also by {\em Fueter elements}, so in our notation, a Fueter element is
an element $\psi = \sum_{c \in P^*} n_c \sigma_c^{-1}$ with $n_c \geq 0$ and $n_c + n_{p-c} = 1$; in particular, $\psi + \jmath \psi = \Norm$.

We can write any $\theta \in I$ as
\begin{eqnarray}
\label{fueterbase} 
 \theta = \sum_{n=1}^{(p-1)/2} \nu_n \psi_n = \sum_{c=1}^{p-1} n_c \sigma_c^{-1}; \quad \nu_n, n_c \in \Z. 
 \end{eqnarray}

Therefore, for each ideal $\eu{C} \subset \Z[ \zeta ]$ and each $\theta \in I$, the ideal $ \eu{C}^{\theta} \subset  \Z[ \zeta ]$ is generated by some $\gamma \in \Z[\zeta ]$, which satisfies $\gamma \cdot \overline{\gamma} = \Norm( \eu{C} )^{\varsigma_{\theta}}$, for an integer $\varsigma_{\theta} \in \Z$, which we call the {\em relative weight} of $\theta$.
 The {\em absolute weight} (or simply, {\em weight},) of  $\theta = \sum_{c \in P^*} n_c \sigma^{-1}_c \in \Z[ G ]$ is $w( \theta ) = \sum_c | n_c |$. We say that $\theta = \sum_{c \in P^*} n_c \sigma_c^{-1}$ is {\em positive}, writing $\theta \in I^+$, if $n_c \in \Z_{\geq 0}$ for all $c \in P^*$. Note that the relative weight of each Fueter element is $1$.

\subsubsection{The Fermat quotient ideal and $J_k \subset I$}
We define the Fermat quotient map $\phi : \Z[ G ] \ra \F_p$ such that
$\zeta^{\theta} = \zeta^{\phi( \theta)}$. Explicitly,
\begin{eqnarray}
\label{phi}
\phi \left(  \sum_{c \in P^*} n_c \sigma_c^{-1} \right) = \sum_{c \in P^*} n_c/c \in \F_p. 
\end{eqnarray}
We identify the value $\phi( \theta ) \in \F_p$ with its natural lift to $\N$,
under the least positive remainder representation of $\F_p$. 
\begin{definition}
\label{FID}
    The {\em Fermat ideal} is $I_0 = I \cap \Ker( \phi )$:
this is the module of all Stickelberger elements $\theta$ such that $\zeta^{\theta} = 1$.
The module $J_k \subset I_0$ is defined by
\begin{eqnarray}
    \label{jk}
    J_k = \{ \theta \in I_0^+ \ : \ \varsigma( \theta ) = 2k, \ k \geq 1 \}.
\end{eqnarray}
It is thus the submodule of $I_0$ consisting of elements that are sums $\theta = \theta_1 + \theta_2$
with $\varsigma( \theta_i ) = k; i = 1, 2$; the $\theta_i$ need not be elements of $I_0$.
\end{definition}
The following fact shows that we always can choose elements $\theta \in I_0^+$
of small relative weight:
\begin{fact}
 \label{i0}
 For $p \geq 5$ there always exists an element $\theta \in I_0^+$ with 
 $\varsigma_{\theta}\geq 2$.
\end{fact}
\begin{proof}
 Let $\phi( \psi_1 ) = a$ and $\phi( \psi_2 ) = b$. 
 If $a \cdot b \equiv 0 \bmod p$, then there exists $j \in \{ 1, 2 \}$ 
 such that $\phi( \psi_j ) = 0$, so $\theta = 2 \psi_j$
 satisfies the claim. Otherwise, let $c \in \N$ be such that
 $ a + b \cdot c \equiv 0 \bmod p$. Since $\zeta^{\phi( \theta ) } = \zeta^{\theta}$, 
 it follows that $\phi( \sigma_c \theta ) \equiv c \phi( \theta ) \bmod p$, 
 and thus $\phi( \psi_1 + \sigma_c \psi_2 ) = 0$. 
 Therefor, $\theta = \psi_1 + \sigma_c \psi_2$ satisfies the claim. 
 Since $I^-$ contains for $p \geq 5$ at least two $\Z$-base elements, the claim follows.
 The same procedure can be applied a fortiori for 
 larger relative weights. 
\end{proof}

\subsubsection{The action of the Stickelberger ideal on characteristic ideals}
\label{sti}
Since $\theta \in I$ annihilates the class group, there is some principal ideal $b(\theta)$ such that
$b = \eu{A}^{\theta}$ and it satisfies $b(\theta) \cdot \bar{b}(\theta) = \Norm( \eu{A} ) = ( z )$.
It is known from the theory of Gauss and Jacobi sums  -- see e.g.
\cite{Wa}, \S 15.1 -- that principal ideals arising from the action 
of the Stickelberger ideal are generated by 
{\em Jacobi numbers}, which are defined as products of Jacobi sums.

Iwasawa proved in \cite{Iw} that Jacobi numbers $\rg{J}$ verify
\footnote{In the classical definition $\tau(\chi) = \sum_{x \in \ZMs{q}} \chi( x ) \xi^x$,
Iwasawa actually proves that $\tau(\chi) \equiv -1 \bmod \eu{P}$, with $\eu{P} \in \Z[ \zeta_p, \xi ]$ 
an ideal above $p$. This led Lang to modify the definition by changing the sign, 
and we use his definition here.}

\begin{eqnarray}
\label{Iw}
 \rg{J} \equiv 1 \bmod ( 1 - \zeta )^2,
\end{eqnarray} 
Since the product of Jacobi numbers by their complex conjugates are rational integers, the above
condition implies that there is a unique Jacobi number that generates the ideal $b(\theta)$, and all
other generators of this ideal, which are rational upon multiplication by their complex conjugates, 
differ from the Jacobi number in $b$ by a root of unity -- a consequence
of the Kronecker unit theorem. 

Let $\beta \in b(\theta)$ be the unique Jacobi number generating the ideal $b(\theta)$. Since 
$(\alpha^{\theta}) = \eu{A}^{p \theta} = b(\theta)^p$, we obtain by using Lemma \ref{aux1}, point 4, 
in combination with \rf{Iw}, that
\begin{eqnarray}
\label{upto}
\alpha^{(1-\jmath) \theta} = \eta' \cdot {\beta}^p, \quad
\eta' = \begin{cases}
         \zeta^{-2x \theta/(x+y)} & \hbox{ for $e = 0$, and}\\
         (-\zeta^{-\theta}) & \hbox{otherwise.}
        \end{cases}
\end{eqnarray}

\subsubsection{The $\beta$-map}
\label{betamap}
With this, we let for $\theta \in I$, the auxiliary number $\beta = \beta( \theta )$ 
be the Jacobi number generating $\eu{A}^{\theta}$.
\begin{lemma}
    \label{betahom}
    The map $\beta : I \ra \K^{\times}$ defines an injective homomorphism of $G$-groups, via 
    \begin{eqnarray}
    \label{bhom}
         \beta( \theta_1 + \theta_2 ) = \beta( \theta_1 ) \cdot \beta( \theta_2 ); \quad \beta( - \theta ) = 1/\beta( \theta ).
    \end{eqnarray}
\end{lemma}
\begin{proof}
    The homomorphism relations in \rf{bhom} can easily be verified from the definition. 
    For $\sigma \in G$, we have by definition 
    \[ \beta( \sigma \theta ) = \beta( \theta )^{\sigma} \in \K^{\times},\]
    so the homomorphism is one of $G$-groups. For injectivity, 
    suppose that $\beta( \theta ) = 1$. Then 
    \begin{eqnarray*}
        \beta( \theta )^q  =  {\alpha'}^{\theta} = \prod_{c=1}^{p-1} \sigma_c^{-1}( \alpha' )^{n_c} = 1.
    \end{eqnarray*}
    By Lemma \ref{aux1}, the ideals $\left( \sigma_c^{-1}( \alpha' ) \right)$ are pairwise 
    coprime, so the exponents $n_c$ must all vanish.
\end{proof}

\subsubsection{Some properties of the module $J_k$}
\label{sjk}
We start by noticing that for $t \in J_k$ we always have ${\alpha'}^t = \alpha^t$ and thus
\begin{eqnarray}
\label{pt}
\beta( t )^q = \alpha^t. 
\end{eqnarray}
This is a consequence of $\zeta^t = 1$. 
We note the following actions $t \in J_k$ on $\lambda$:
\begin{eqnarray}
\label{pm}
    ( 1 - \zeta )^{t( 1+ \jmath)}  =  p^{2k}; \quad ( 1 - \zeta )^{t( 1 - \jmath)} = (-\zeta)^t = (-1)^{k(p-1)} \cdot \zeta^t = 1, \quad \hbox{hence} \
    \lambda^t = s(t)  p^k, \ s( t ) \in \{ -1, 1 \}. 
\end{eqnarray}

\begin{definition}
    \label{alredef}
    In view of \rf{pt} and \rf{pm}, we shall from now on work only with $t \in J_k$
    and redefine $\alpha = (y + \zeta x )$ also for the case $e = 1$. This allows a unified
    treatment of both cases of \rf{sfct}, when $( K, p ) = 1$ and when $p | K$.
\end{definition}
In the case $e = 1$, the first line in \rf{pm} implies that
$\beta^q = s(t) \alpha^{\theta}/p^{k}$ and consequently
\begin{eqnarray}
\label{betasm1}
\beta^{\sigma-1} = \alpha^{\theta( \sigma-1 )} , \quad 
 \forall \sigma \in G \setminus \{1\}.
\end{eqnarray}
showing that we the denominator vanishes in the case $e = 1$, when acting with $\sigma - 1$ on $\beta( t )$,
for arbitrary $\sigma \in G \setminus \{ 1 \}$ and $t \in J_k$.
The following result induces the key factorizations of $K$, mentioned in the introduction.
\begin{lemma}
    \label{psipol}
    Let $\Psi_t( x, y ) = \alpha^t$ and $K = x y (x+y )$. Then 
    $K | \left( \Psi_t (x, y ) - \sigma(\Psi_t (x, y ))\right)$ for all $\sigma \in G  \setminus \{1\}$.
\end{lemma}
\begin{proof}
    We have the following congruences:
    \begin{eqnarray*}
        \Psi_t( x, y ) & \equiv & \begin{cases}
            y^{(p-1) k} \bmod x, \\
            x^{(p-1) k} \bmod y, \\
            \lambda^t y^{(p-1)k} = s(t) p^k y^{(p-1) k} \bmod x+y; \ s( t ) \in \{ -1 , 1 \}
        \end{cases}
    \end{eqnarray*}
    The right side of the congruence being rational in all cases, the claim follows.
\end{proof}

\section{The Strong Fermat-Catalan Equation}
We assume that \rf{sfct} has a non trivial solution, let 
$\alpha = y + \zeta x$ for both values as $e = 1$, as mentioned above.
Let $t \in J_k$. Let $\beta( t )$ be the Jacobi sum generating $\eu{A}^t$.

We define the following:
\begin{definition}
\label{datadef}
Let $t \in J_k$ and $w \in \{ x, y, x+y\}$. For each pair $(t, w)$ and $\sigma \in G \setminus \{ 1 \}$, let:
\label{data1}
 \begin{eqnarray*}
     \Delta( \sigma ) & = & \Delta_{t}( \sigma )  =  \Psi_t - \sigma \Psi_t, \\
     \delta_c(\sigma ) & = & \beta - \xi^c \sigma(\beta), \quad c \in Q, \\
      w' & = & \frac{w}{q^{v_q( w )}}, \\    
     \eu{D}_c(\sigma ) & = & ( w', \delta_c( t, \sigma ) ), \quad c \in Q, \\
     \rg{D} & = & \bigcap_{\sigma \in G} \eu{D}_0( \sigma ).
 \end{eqnarray*}    
\end{definition}
The items in the above definition depend on a set of variables: $t, w, \sigma, c$. In the 
sequel, we shall always keep the reference on $c$ explicit, while the further variables
will only mentioned when their choice changes or is not obvious in the contest.

As an immediate consequence of the definitions and of Lemma \ref{aux1}, we have
\begin{lemma}
    \label{cop}
    Let $t \in J_k; w \in \{ x, y, x+y\}$ and $\sigma \in G \setminus \{ 1 \}$ be fixed.
    Then
    \begin{itemize}
    \item[ 1. ] For distinct $a, b \in Q$, the ideals $\eu{D}_a, \eu{D}_b$ are pairwise coprime.
    \item[ 2. ] The product  $\prod_{a \in Q} \eu{D}_a = ( w' )$ and $( \eu{D}_a, q ) = 1$. 
    \item[ 3. ] The ideal $\rg{D}$ is $G$-invariant, so $\rg{D} = (D)$ for some $D \in \N$.
    Moreover, $\beta( t ) \equiv \sigma( \beta ( t )) \equiv z^k \bmod \rg{D}$
    for all $\sigma \in G$.
    \item[ 4. ] Let $\eu{R} \subset \Z[ \zeta ]$ be a prime that divides $w'$. If $\eu{R} \nmid \rg{D}$, then there is a $\sigma \in G$ and $a \in Q^*$ such that $\eu{R} | \eu{D}_a( \sigma )$.
    \end{itemize} 
\end{lemma}
\begin{proof}
    Defining the ideal $D( a, b ) = ( \delta_a, \delta_b )$, we see for $a, b > 0$ that
    \begin{eqnarray*}
        (\xi^b - \xi^c ) \sigma( \beta ) &\in& D( a, b ), \quad \hbox{ and } \ 
        \ol{(\xi^b - \xi^c )} \beta \in D( a, b ) \ \quad \hbox{hence} \\
        D( a, b ) & \mid & ( \beta, \sigma( \beta )) \cdot (\lambda' ).
    \end{eqnarray*}
    Now $\beta | z^{\varsigma( t )}$ and since $( z, w ) = 1$ it follows that 
    $( \eu{D}_a, \eu{D}_b ) | ( \lambda' )$.
    By definition however, $( \eu{D}_a, q ) = 1$ so we conclude 
    that $( \eu{D}_a, \eu{D}_b ) = 1$.
    A similar argument holds when one of $(a, b )$ is $0$ -- one may assume $a = 0$ 
    and build appropriate differences to show that $( \eu{D}_0, \eu{D}_b ) = 1$ 
    in this case too; we leave the details to the reader.

    We have shown that $\Psi_t(x, y ) \equiv \sigma( \Psi_t( x, y ) \bmod w$; 
    since $\prod_{a \in Q} \delta_a = 
    \Psi_t( x.y ) - \sigma(\Psi_t( x, y )$, while $( \eu{D}_a, q ) = 1$, it follows that
    $\prod_{a \in Q} \eu{D}_a = (w')$, which confirms 2.

    From the definition, if a prime $\eu{R} | \rg{D}$, then 
    $\beta \equiv \sigma( \beta ) \bmod \eu{R}$ for 
    all $\sigma \in G$. Thus $\beta \equiv C = \frac{\Tr( \beta ) }{p-1} \bmod \eu{R}$, 
    with $C \in \Q$.
    Since $\beta \equiv \bar{\beta} \bmod \eu{R}$, it also follows by multiplying the 
    two congruences, that $z^{2k} \equiv C^2 \bmod \eu{R}$, hence 
    $C \equiv s z^k \bmod \eu{R}$, for $s \in \{-1,1\}$. 
    We claim that the sign of $C$ is in fact positive. For some $w \in \{ x, y, x + y, x-y\}$, 
    and $\eu{R} \ | \ \rg{D} \ | \ w'$, we have, by raising to the odd power $q$,
    \begin{eqnarray*}
        \beta( \sigma, w ) & \equiv & s z^k \bmod \eu{R} \quad \Rightarrow \quad \alpha^{\sigma t} \equiv s z^{qk} \bmod \eu{R}.
    \end{eqnarray*}
    Using Lemma \ref{psipol}, we note that $k \cdot \Norm \in J_k$ 
    and the proof of the Lemma also applies to 
    $z^{ q k} = \Norm( \alpha^{k \sigma })$. Consequently, $ \alpha^{\sigma t} \equiv z^{q k} \bmod K$, and thus $s = 1$, since $q$ is odd. This confirms the claim.

    By acting with $G$ on the congruence we derived, we conclude that:
    \begin{eqnarray}
        \label{valD}
        \sigma( \beta ) \equiv z^k \bmod \tau( \eu{R} ), \quad \forall \sigma, \tau  \in G.
    \end{eqnarray}
    The congruences thus hold modulo $r$, the rational prime below $\eu{R}$. 
    This holds for all $\eu{R} | \rg{D}$, so $\rg{D}$ is $G$-invariant. This confirms 3.
    
    Moreover, any prime $\eu{R} | ( w' )$ which is coprime to $\rg{D}$
    will divide, by the way $\rg{D}$ was defined, 
    $\eu{D}_a( \sigma )$ for some $a \in P^*$ and $\sigma \in G \setminus \{ 1 \}$. 
    This confirms point 4.
\end{proof}
The above Lemma indicates that the $\eu{D}_a$ induce uncommonly high factorizations of
the numbers $w \in F =\{ x,y, x+y, x-y\}$. This observation is not particularly new,
and in itself it cannot bring final insights about the equations under investigation.
However, the definition of the ideal $\rg{D}$ leads further, and we deduce from the above:
\begin{lemma}
    \label{gcd0}
    Under the notations of Definition \ref{datadef}, we have $D = 1$.
\end{lemma}
\begin{proof}
    By Fact \ref{i0}, we can choose $t \in J_k$ for any $p > 3$ and $k \leq \frac{p-1}{2}$.
    We may write $t = t_1 + t_2$, with $t_i \in I( k ); i = 1, 2$. Then
    $t = t_1 + t_2 = t_1 + k \Norm - \bar{t}_2$.
    
    For any $w'$ in Definition \ref{datadef}, and any $\sigma \in G \setminus \{ 1 \}$, we know 
    from \rf{valD} that $\sigma( \beta( t )) \equiv z^k \bmod D$. 
    There is thus a $\chi \in \Z[ \zeta ]$ such that
    \[ \beta( t ) -  z^k = D \cdot \chi. \]
    Write now $\beta( t ) = \beta( t_1 ) \cdot \beta( t_2 )$ and 
    $z^k = \beta( t_2 )^{1+\jmath}$. We then have
    \[ \beta( t_1 ) \cdot \left( \beta( t_2) - \bar{\beta}( t_1 )\right) = D \cdot \chi. \]
    Since $D | w$ and $\beta | z^{2k}$, it follows that $( \beta( t_1 ), D ) = 1$, 
    so we conclude that $\beta( t_1 ) | \chi$. We repeat the same argument, by 
    interchanging the roles of $t_1$ and $t_2$, finding that $\beta( t_2 ) | \chi$ 
    too. There is thus a $\chi' \in \Z[ \zeta ]$, such that 
    $\chi = \beta( t ) \chi'$. The initial identity becomes
    \[ \beta( t ) \cdot ( 1 - D \chi' ) = z^k .\]
    Multiplying by complex conjugates in the same identity, we
    conclude that $( 1 - D \chi' )^{1+\jmath} = 1$, and thus 
    $\mu := 1 - D \chi' \in \lan \pm \zeta \ran$, by the Kronecker Unit Theorem. 
    So $\beta( t ) = \bar{\mu} z^k$ and $\mu^{2p} = 1$.
    But then $( \bar{\mu} z^k )^{2 p q} \in \Z$, while $\beta( t )^{2p q} = \alpha^{2t p} \not \in \Z$,
    as follows from Lemma \ref{aux1}. 
    The assumption $D \neq 1$ is thus untenable; of course, if $D = 1$, 
    the original congruences modulo $D$ are void, so there is no contradiction.
\end{proof}

The Lemma implies that for any $t \in J_k$, every rational factor of $w'$
will have a non trivial factor in some ideal $\eu{D}_a( \tau )$, 
for some $\tau \in G \setminus \{ 1 \}$. This fact, together with the result of Lemma \ref{gcd0}
allows us to show that these factors split completely in $\K'$:
\begin{lemma}
    \label{gcd1}
    Let $w' \in F' = \{ x', y', (x+y)' \}$, as defined in Definition \ref{datadef}.
    Then every prime $r | w'$ is totally split in $\K'$.
\end{lemma}
\begin{proof}
Let $t = t_1 + t_2 \in J_1$ and consider a prime $r | w'$. Since $D = 1$, there is a
$\sigma \in G$ such that $\left( r, \frac{(w')}{\eu{D}_0( \sigma )}\right) \neq (1)$.
This implies that $\frac{r}{(r, \eu{D}_0( \sigma )}$ is a non trivial 
product of primes of $\K$ which split completely in $\KL/\K$. 
If $\rg{r} \subset \id{O}( \KL )$ is a prime above $\eu{R}$, and 
$D_r \subset \Gal( \KL/\Q )$ is its decomposition group, it follows that 
$D_r \cap H' = \{ 1 \}$ and thus $\KL^{D_r} \supset \K'$. Consequently, the rational
prime $r$ below $\eu{R}$ splits completely in $\K'/\Q$, which completes the proof.
\end{proof}
We thus conclude that
\begin{proposition}
    \label{pmain}
    Let $K' = \prod_{w | xy (x+y )} w'$.
    The primes $r | K'$ are all totally split in $\K'$, so $r \equiv 1 \bmod q$.
    Moreover, if $r | K$ then $r | K'$ or $r = q$.
\end{proposition}
\begin{proof}
    Let $F = \{ x, y, x+y\}$ and $' : F \ra \Z$ the map 
    $v \mapsto v/q^{v_q( v )}$ that 
    sends $w \in F$ to $w'$ as defined in Lemma \ref{cop}. We have 
    $K = \prod_{w \in F} w$, so $K$ and $K'$ differ by a power of $q$. 
    By Lemma \ref{gcd1}, the primes $r | K'$ are totally split in $\K'$, 
    thus verifying $r \equiv 1 \bmod q$. The second claim follows from
    the definition of $K'$.
\end{proof}
As a consequence:
\begin{corollary}
    \label{thm}
    Equation \rf{sfct} has no solutions for $p > 3$ and Theorem \ref{main} is true.
\end{corollary}
\begin{proof}
    At least one $v \in F$ must be even: indeed, if $x$ and $y$ have the same parity, 
    then $x + y \equiv 0 \bmod 2$. Otherwise, one of $x$ or $y$ must be even, 
    thus confirming the claim. Since $q$ is odd, $v$ and $v'$ have 
    the same parity for all $v \in F$; it remains that $2 | K'$. By Proposition
    \ref{pmain}, all primes $r | K$
    are of the form $r = m q + 1$ or $r = q$, so this should hold also for 
    $r = 2$. We reach a contradiction, which confirms Theorem \ref{main}.
\end{proof}

\vspace*{0.5cm}
\textbf{Acknowledgements}
This work is built on shoulders of giants, mainly on the theories and results due to Kummer, Hensel and Chevalley. 
One cannot however under-appreciate the contribution of numerous later authors, whose work with these methods made a combination of ideas which might have seemed remote in earlier times, become possible, in a natural way. The results of
Furtw\"angler from the early \nth{20} century appear particular as a forecast of
the present result.

I am deeply indebted to the people who contributed to this recent surprising development. We are grateful to  Ivan Fesenko and the students H. Chen, D. Lehmann, C. Liu and J. Reichardt for their active interest and questions which helped improve the presentation. 

It is not possible to list, or even
recall, the sporadic yet important discussions and challenges which made that, despite the decision to stay away from ternary equations, the interest in the question was warmed up periodically, so it never completely disappeared, over more than two decades. The  support of the colleagues at the Mathematical Institute of the University of G\"ottingen will be mentioned, collectively -- to all, our thanks.

Most of all, I am grateful to my family in small -- Theres and Seraina -- and at large, unlistable, who stood beside in support during short and 
long times.

\end{document}